\input amstex
\documentstyle{amsppt}
\magnification=\magstep1                        
\hsize6.5truein\vsize8.9truein                  
\NoRunningHeads
\loadeusm

\magnification=\magstep1                        
\hsize6.5truein\vsize8.9truein                  
\NoRunningHeads
\loadeusm

\document
\topmatter

\title
Extensions of the  Bloch-P\'olya theorem on the number of real zeros of polynomials (II)
\endtitle

\author
Tam\'as Erd\'elyi
\endauthor

\address Department of Mathematics, Texas A\&M University,
College Station, Texas 77843 \endaddress

\email terdelyi\@tamu.edu
\endemail

\thanks {{\it 2020 Mathematics Subject Classifications.} 26C10, 11C08}
\endthanks

\keywords location of zeros, polynomials with constrained coefficients
\endkeywords

\date October 13, 2024 \enddate

\abstract
We prove that there is an absolute constant $c > 0$ such that for every 
$$a_0,a_1, \ldots,a_n \in [1,M]\,, \qquad 1 \leq M \leq \frac 14 \exp \left( \frac n9 \right),$$
there are 
$$b_0,b_1,\ldots,b_n \in \{-1,0,1\}$$ 
such that the polynomial $P$ of the form $\displaystyle{P(z) = \sum_{j=0}^n{b_ja_jz^j}}$ 
has at least $\displaystyle{c \left( \frac{n}{\log(4M)} \right)^{1/2}-1}$
distinct sign changes in $I_a := (1-2a,1-a)$, where 
$\displaystyle{a := \left( \frac{\log(4M)}{n} \right)^{1/2} \leq 1/3}$. 
This improves and extends earlier results of Bloch and P\'olya and Erd\'elyi
and, as a special case, recaptures a special case of a more general recent 
result of Jacob and Nazarov.
\endabstract

\endtopmatter

\head 1. Introduction \endhead

In this note $c$ and $c_j$ will denote suitable positive absolute constants.
Let $\Cal{F}_n$ denote the set of polynomials of degree at most $n$
with coefficients from $\{-1,0,1\}$. Let $\Cal{L}_n$ denote the set of
polynomials of degree $n$ with coefficients from $\{-1,1\}$.
Let $\Cal{A}_n$ denote the set of polynomials of degree at most $n$
with coefficients from $\{0,1\}$. Updating what the authors write in 
[BEK99] and [Er08] we have the following. The study of the location 
of zeros of polynomials from the above classes begins with Bloch and P\'olya 
[BP32]. They prove that the average number of real zeros of a polynomial 
$0 \not\equiv P \in \Cal{F}_n$ is at most $cn^{1/2}$. They also prove 
that a polynomial $0 \not \equiv P$ from $\Cal{F}_n$ cannot have more than 
$$\frac{cn\log\log n}{\log n}\,, \qquad  n \geq 3\,,$$ 
real zeros. This quite weak result appears to be the first on this subject. 
Schur [Sch33] and by different methods Szeg\H o [Sze34] and 
Erd\H os and Tur\'an [ET50] improve this upper bound to $c(n\log n)^{1/2}$, 
see also [BE95]. (Their results are more general, but in this specialization 
not sharp.) Our Theorem 4.1 in [BEK99] gives the right upper bound, $cn^{1/2}$, 
for the number of real zeros of polynomials 
from a much larger class, namely for all polynomials $P$ of the form
$$P(x) = \sum_{j=0}^n{a_jx^j}\,, \quad |a_0| = |a_n| = 1\,, \quad |a_j| \leq 1\,, 
\quad a_j \in {\Bbb C}\,.$$
Schur [Sch33] claims that Schmidt gives a version of part of this theorem.
However, it does not appear in the reference he gives, namely [Sch32], and we
have not been able to trace it to any other source. In [BE97]
we prove that $\eta n^{1/2}$ is an upper bound for the number of zeros of 
polynomials $P$ of the form 
$$P(x) = \sum_{j=0}^n{a_jx^j}\,, \quad |a_0| = 1\,, \quad |a_j| \leq 1\,, \quad a_j \in {\Bbb C}\,, \tag 1.1$$
inside any polygon with vertices on the unit circle, where the constant $\eta$ 
depends only on the polygon. Bloch and P\'olya [BP32] also proved that there are 
polynomials $P \in \Cal{F}_n$ with at least  
$$\frac{cn^{1/4}}{(\log n)^{1/2}}\,, \qquad n \geq 2\,,$$
sign changes on the real line ${\Bbb R}$. 
Schur [Sch33] claims they do it for polynomials with coefficients only from $\{-1,1\}$, 
but this appears to be incorrect. 

In a seminal paper Littlewood and Offord [LO39] prove that the number of real
roots of a polynomial $P \in \Cal{L}_n$, on average, lies between
$$\frac{c_1\log n}{\log\log\log n} \qquad \text{and} \qquad c_2(\log n)^2, \qquad n \geq 3\,,$$
and it is proved by Boyd [Bo97] that every $P \in \Cal{L}_n$ has at most
$\displaystyle{\frac{(1+o(1))(\log n)^2}{\log\log n}}$ zeros at $1$ (in the sense 
of multiplicities). 

Kac [Ka48] proves that the expected number of real roots of a polynomial of degree $n$
with independent and identically uniformly distributed coefficients on $[-1,1]$ is 
asymptotically $(2/\pi)\log n$.

In [Er08] we prove that there are absolute constants
$c_1 > 0$ and $c_2 > 0$ such that for every
$$a_0,a_1, \ldots,a_n \in [1,M]\,, \qquad 1 \leq M \leq \exp(c_2n^{1/4})\,,$$
there are $b_0,b_1,\ldots,b_n \in \{-1,0,1\}$ such that
$\displaystyle{P(z) = \sum_{j=0}^n{b_ja_jz^j}}$ has at least $c_1n^{1/4}$
distinct sign changes in $(0,1)$. This improves and extends the earlier result
of Bloch and P\'olya [BP32].

It is proved in [BEK99] (see Theorem 2.1) that every polynomial $P$ of the form 
$$P(x) = \sum_{j=0}^n{a_jx^j}\,, \quad |a_0| = 1\,, \quad |a_j| \leq M\,,\quad a_j \in {\Bbb C}\,, 
\quad M \geq 1\,, \tag 1.2$$ 
has at most $c_1(n(1+\log M))^{1/2}$ zeros (in the sense of multiplicities) at $1$.
It is also proved in [BEK99] (see Theorem 4.1) that every polynomial $P$ of the form (1.1) 
has at most $c_2n^{1/2}$ zeros in $[-1,1]$. Modifying the proof of Lemma 5.6 in [BEK99] 
(see also the proof of Lemma 1 on page 55 of [Bo02]) Jacob and Nazarov prove in the 
Appendix of [JN24] that every polynomial $P$ of the form (1.2) has at most $c_3n^{1/2}(1+\log M)$ 
zeros in $[-1,1]$. 

Jacob and Nazarov [JN24] prove a general theorem including that
there are polynomials $P \in \Cal{L}_n$ having at least $c_1n^{1/2}$ distinct zeros
in $[0,1]$; and there are polynomials $0 \not \equiv P \in \Cal{A}_n$ for which the
polynomials $Q$ defined by $Q(x) := P(-x)$ have at least $c_2n^{1/2}$ distinct zeros
in $[0,1]$.

\head 2. New Result \endhead

\proclaim{Theorem 2.1} 
There is an absolute constants $c > 0$ such that for every
$$a_0,a_1, \ldots,a_n \in [1,M]\,, \qquad 1 \leq M \leq \frac 14 \exp \left( \frac n9 \right), \tag 2.1$$
there are 
$$b_0,b_1,\ldots,b_n \in \{-1,0,1\}$$ 
such that the polynomial $P$ of the form $\displaystyle{P(z) = \sum_{j=0}^n{b_ja_jz^j}}$
has at least $\displaystyle{c \left( \frac{n}{\log(4M)} \right)^{1/2}-1}$ 
distinct sign changes in $I_a := (1-2a,1-a)$, where 
$\displaystyle{a := \left( \frac{\log(4M)}{n} \right)^{1/2} \leq 1/3}$.
\endproclaim

This improves and extends earlier results of Bloch and P\'olya [BP32] and Erd\'elyi [Er08],
and, as a special case, recaptures the special case of a more general recent result of 
Jacob and Nazarov [JN24] mentioned above.

\head 3. Lemmas \endhead

We use the notation $\log^{-}y := \max\{0, -\log y\}$ and $\log^{+}y := \max\{0, \log y\}$ for $y \in (0,\infty)$.
The geometric mean $L_0([a,b])$ of the absolute value of a real-valued continuous function $P$ on a finite interval 
$[a,b] \subset {\Bbb R}$ is
$$\|P\|_{L_0([a,b])} := \exp \left( \frac{1}{b-a}\int_a^b{\log|P(x)| \,dx} \right)\,.$$
In [JN24] Jacob and Nazarov showed the following.      

\proclaim{Lemma 3.1}
Let $a \in (0,1/3]$ and $I_a := (1-2a,1-a)$. Then
$$\int_{I_a}{\log^{-}|P(x)| \, dx} \leq c_1 \log(4M)$$ 
for every polynomial $P$ of the form  
$$P(x) = \sum_{j=0}^{n}{b_ja_jx^j}\,, \quad a_j \in [1,M]\,, \quad b_j \in \{-1,0,1\}\,, \quad b_0 \in \{-1,1\}\,. \tag 3.1$$ 
Moreover $c_1 := 9\pi/\sqrt{2}$ is a suitable choice.
\endproclaim

As a consequence we have the following.

\proclaim{Lemma 3.2} Let $a \in (0,1/3]$ and $I_a := (1-2a,1-a)$. Then
$$\|P\|_{L_0(I_a)} \geq \exp \left( \frac{-c_1\log(4M)}{a} \right)$$
for every polynomial $P$ of the form (3.1).
\endproclaim 

We will need a simple consequence below of Lemma 3.2 in the proof of Theorem 2.1.

\proclaim{Lemma 3.3}  
Let $a \in (0,1/3]$, $I_a := (1-2a,1-a)$, and $0 \leq k \leq n$. Then
$$\|P\|_{L_0(I_a)} \geq (1-2a)^k \exp\left( \frac{-c_1\log(4M)}{a} \right)$$
for every polynomial $P$ of the form 
$$P(x) = \sum_{j=k}^{n}{b_ja_jx^j}\,, \quad a_j \in [1,M]\,, \quad b_j \in \{-1,0,1\}\,, 
\quad b_k \in \{-1,1\} \,. \tag 3.2$$
\endproclaim 

\head 4. Proofs \endhead

\demo{Proof of Lemma 3.2}
Using Lemma 3.1 we have 
$$\split \|P\|_{L_0(I_a)} = &  
\exp \left( \frac{1}{a} \int_{I_a}{\log^{+}|P(x)| \, dx} - \frac{1}{a} \int_{I_a}{\log^{-}|P(x)| \, dx} \right) \cr 
\geq & \exp \left(\frac{-1}{a} \int_{I_a}{\log^{-}|P(x)| \, dx} \right) \cr 
\geq & \exp \left( \frac{-c_1\log(4M)}{a}\right)\,. \cr 
\endsplit$$ 
\qed \enddemo 

\demo{Proof of Theorem 2.1}
We modify the outline given by Nazarov in MathOverflow [MO24] where $a_j=1$ for each $j=0,1,\ldots,n$, is assumed.
Let $a \in (0,1/3]$ and $I_a := (1-2a,1-a)$.
Associated with a polynomial $P$ we define the polynomial $\widetilde{P}(x) := P(1-a-ax)$.
The function $f(x) := 1-a-ax$ maps the interval $(0,1)$ onto the interval $(1-2a,1-a)$, hence 
the polynomial $P$ of the form 
$$P(x) = \sum_{j=0}^n{d_ja_j x^j}, \qquad a_j \in [1,M]\,, \quad d_j \in \{0,1\}\,,  \tag 4.1$$ 
satisfies
$$0 \leq \widetilde{P}(x) \leq (n+1)M\,, \qquad x \in [0,1]\,.$$ 
For every fixed $a_0,a_1,\ldots,a_n \in [1,M]$ there are $2^{n+1}$ different polynomials $P$ of the form (4.1).
Let $1 \leq m \leq n$ be integers. We study the $m$-tuples 
$$\left(\int_0^1{\widetilde{P}(x)x^j\,dx} \right)_{j=0}^{m-1}$$
which are in the $m$-dimensional cube $[0,(n+1)M]^m$.  
We divide the cube $[0,(n+1)M]^m$ into $N=L^m$ equal subcubes, so the side length of the subcubes is 
$h := (n+1)M/L$. We choose $L := [2^{n/m}]$, where $[x]$ denotes the integer part of a nonnegative real number $x$. 
Then $N = L^m \leq 2^n < 2^{n+1}$ and $L \geq 2^{n/m}/2$,  
hence the side length $h$ of the subcubes satisfies $h = (n+1)M/L \leq 2(n+1)M2^{-n/m}$.
By the pigeonhole principle there are two different polynomials $P_1$ and $P_2$ of the form (4.1), 
that is, 
$$P_1(x) = \sum_{j=0}^n{d_ja_jx^j}, \qquad d_j \in \{0,1\}\,,$$
and
$$P_2(x) = \sum_{j=0}^n{d_j^*a_jx^j}, \qquad d_j^* \in \{0,1\}\,,$$ 
such that for some $0 \leq k \leq n$ the polynomial $0 \not\equiv P := P_1-P_2$ is of the form (3.2) with 
$$b_j := d_j-d_j^* \in \{-1,0,1\}\,, \qquad j=k,k+1,\ldots,n \,,$$  
and 
$$\left|\int_0^1{\widetilde{P}(x)x^j\,dx}\right| \leq 2(n+1)M2^{-n/m}\,, \qquad j=0,1,\ldots,m-1\,. \tag 4.2$$ 
For $0 \leq k \leq n$ we have 
$$(1-2a)^k \geq (1-2a)^n \geq \exp(-4an)\,, \quad a \in (0,1/3]\,.$$ 
Combining this with Lemma 3.3 gives
$$\|P\|_{L_0(I_a)} \geq \exp \left( \frac{-c_1\log(4M)}{a}-4an \right)\,.$$
Now we choose $\displaystyle{a := \left( \frac{\log(4M)}{n} \right)^{1/2} \leq 1/3}$ to get 
$$\|\widetilde{P}\|_{L_0([0,1])} = \|P\|_{L_0(I_a)} \geq \exp \left( -(c_1+4)(\log(4M))^{1/2}n^{1/2} \right)\,. \tag 4.3$$ 
Assume that $P$ has exactly $m-1 \geq 0$ sign changes in $I_a$. Then $\widetilde{P}$ has 
exactly $m-1$ sign changes in $(0,1)$. We denote the sign changes of $\widetilde{P}$ in $(0,1)$ 
by $x_1,x_2,\ldots,x_{m-1}$. We define 
$$Q(x) := \prod_{j=1}^{m-1}{(x-x_j)}$$ 
(the empty product is defined to be $1$). 
Then $\widetilde{P}Q$ does not change sign in $(0,1)$. Observe that   
$$\|x-x_j\|_{L_0([0,1])} = \exp \left( \int_0^1{\log|x-x_j|\, dx} \right) \geq 
\exp \left( \int_0^1{\log x\, dx} \right) = e^{-1}$$
for each $j=1,2,\ldots,m-1$, hence
$$\|Q\|_{L_0([0,1])} \geq e^{-(m-1)} > e^{-m} \,. \tag 4.4$$ 
Combining (4.3), (4.4), the multiplicative property of the geometric mean, the comparison between 
the $L_0([0,1])$ and $L_1([0,1])$ means, the fact that $\widetilde{P}Q$ does not change sign in $(0,1)$, 
(4.2), and the fact that the sum of the absolute values of the coefficients of $Q$ is at most $2^m$, we obtain 
$$\exp \left( -(c_1+4)(\log(4M))^{1/2}n^{1/2} - m \right) \leq \|\widetilde{P}Q\|_{L_0([0,1])} 
\leq \|\widetilde{P}Q\|_{L_1([0,1])}$$ 
$$= \left| \int_0^1{(\widetilde{P}Q)(x)\, dx} \right| \leq 2^m 2(n+1)M 2^{-n/m}.$$
Taking the log of both sides above, we obtain
$$-(c_1+4)(\log(4M))^{1/2}n^{1/2} - m \leq m \log 2 + \log(2n+2) + \log M - \frac nm \log 2\,,$$
that is, 
$$\frac nm \log 2 \leq (\log 2 + 1)m + \log(2n+2) + \log M + (c_1+4)(\log(4M))^{1/2}n^{1/2}\,. \tag 4.5$$
Observe that assumption (2.1) implies 
$$\log M < (\log(4M))^{1/2}(\log(4M))^{1/2} \leq (\log(4M))^{1/2}(1/3)\,n^{1/2}\,,$$
hence it follows from (4.5) that 
$$\frac nm \log 2 \leq (\log 2 + 1)m + \log(2n+2) + (c_1 + 4 + 1/3)(\log(4M))^{1/2}n^{1/2}\,. \tag 4.6$$
Finally, assuming that $\displaystyle{1 \leq m \leq c_2 \left( \frac{n}{\log(4M)} \right)^{1/2}}$ 
in (4.6) with a sufficiently small absolute constant $c_2 > 0$, we get a contradiction.
So $\displaystyle{m \geq c_2 \left( \frac{n}{\log(4M)} \right)^{1/2}}$ with a sufficiently small 
absolute constant $c_2 > 0$. To find such an explicit absolute constant $c_2 > 0$ recall that 
$c_1 := 9\pi/\sqrt{2}$ is a suitable choice. 
\qed \enddemo

\head 5. An open problem \endhead

\proclaim{Problem 5.1} What is the largest multiplicity of the zero a polynomial 
$0 \not\equiv P \in \Cal{F}_n$ can have at $1$? 
\endproclaim

As it is shown by Theorem 4.1 in [BEK99], polynomials
$0 \not\equiv P \in \Cal{F}_n$ cannot have more than $c_1n^{1/2}$ zeros at $1$. By using the pigeonhole principle 
it is shown by Theorem 2.7 in [BEK99] that there exist polynomials $0 \not\equiv P \in \Cal{F}_n$ with at least 
$c_2(n/\log n)^{1/2}$ zeros at $1$. Closing the gap between $c_1n^{1/2}$ and $c_2(n/\log n)^{1/2}$ looks difficult.   

\head 6. Acknowledgment \endhead 
The author thanks Fedor Nazarov for the discussions on the topic and for presenting 
the main idea of the proof of Theorem 2.1 on MathOverflow.


\Refs \widestnumber\key{ABCD2}

\medskip

\ref \no Bo02
\by P. Borwein 
\book Computational Excursions in Analysis and Number Theory
\yr 2002 \publ Spring-Verlag \publaddr New York 
\endref

\medskip

\ref \no BP32 
\by A. Bloch and G. P\'olya
\paper On the roots of certain algebraic equations
\jour Proc. London Math. Soc. 
\yr 1932 \vol 33 \pages 102--114
\endref

\medskip

\ref \no BE95
\by  P. Borwein and T. Erd\'elyi \book
Polynomials and Polynomial Inequalities
\yr 1995 \publ Springer-Verlag \publaddr New York
\endref

\medskip

\ref \no BE97
\by P. Borwein and T. Erd\'elyi
\paper On the zeros of polynomials with restricted coefficients
\jour Illinois J. Math. \vol 41 \yr 1997 \pages 667--675
\endref

\medskip

\ref \no BEK99
\by P. Borwein, T. Erd\'elyi, and G. K\'os
\paper Littlewood-type problems on $[0,1]$.
\jour Proc. London Math. Soc. \vol 79 \yr 1999  \pages 22--46
\endref

\medskip

\ref \no Bo97
\by D.W. Boyd
\paper On a problem of Byrne's concerning polynomials with restricted coefficients
\jour Math. Comp. \vol 66 \yr 1997 \pages 1697--1703
\endref

\medskip

\ref \no Er08
\by T. Erd\'elyi
\paper Extensions of the Bloch-P\'olya theorem on the number of real zeros of polynomials
\jour Journal de th\'eorie des nombres de Bordeaux
\yr 2008 \vol 20, no. 2 \pages 281--287
\endref

\medskip

\ref \no ET50 
\by P. Erd\H os and P. Tur\'an 
\paper On the distribution of roots of polynomials 
\jour Ann. Math. \vol 57 \yr 1950 \pages 105--119 
\endref

\medskip

\ref \no JN24
\by M. Jacob and F. Nazarov
\paper The Newman algorithm for constructing polynomials with restricted coefficients and many real roots
\jour Rev. Uni\'on Mat. Argent. https://inmabb.criba.edu.ar/revuma/pdf/earlyview/4883.pdf  
\yr 2024 \paperinfo https://arxiv.org/pdf/2404.07971 \toappear 
\endref

\medskip

\ref \no Ka48
\by M. Kac
\paper On the average number of real roots of a random algebraic equation (II)
\jour Proc. London Math. Soc. \vol 50 \yr 1948 \pages 390--408
\endref

\medskip

\ref \no LO39
\paper On the number of real roots of a random algebraic equation. II
\by J.E. Littlewood and A.C. Offord
\jour Proc. Cam. Phil. Soc. \vol 35 \year 1939 \pages 133--148
\endref

\medskip

\ref \no MO24
\jour https://mathoverflow.net/questions/461631/number-of-real-roots-of-0-1-polynomial
\endref

\medskip

\ref \no Sch32
\by E. Schmidt
\paper \"Uber algebraische Gleichungen vom P\'olya-Bloch-Typos
\jour Sitz. Preuss. Akad. Wiss., Phys.-Math. Kl. \pages 321 \yr 1932 \endref

\medskip

\ref \no Sch33
\by I. Schur
\paper Untersuchungen \"uber algebraische Gleichungen.
\jour Sitz. Preuss. Akad. Wiss., Phys.-Math. Kl. \pages 403--428 \yr 1933
\endref

\medskip

\ref \no Sze34
\by G. Szeg\H o
\paper Bemerkungen zu einem Satz von E. Schmidt \"uber algebraische Gleichungen.
\jour Sitz. Preuss. Akad. Wiss., Phys.-Math. Kl. \pages 86--98 \yr 1934
\endref

\endRefs
\enddocument